\documentclass[12pt,reqno]{amsart}
\usepackage{amssymb}
\usepackage{amsmath,amsthm}

\newcommand{\be}{\mbox{{\bf E}}}
\newcommand{\bp}{\mbox{{\bf P}}}
\newcommand{\bs}{\mathbf{s}}
\newcommand{\bsp}{\mathbf{s'}}
\newcommand{\bx}{\mathbf{x}}
\newcommand{\bxp}{\mathbf{x'}}

\newcommand{\ot}{[0,t]}

\newcommand{\1}{{\bf 1}}

\newcommand{\D}{\mathbb D}
\newcommand{\R}{\mathbb R}

\newcommand{\Z}{\mathbb Z}
\newcommand{\E}{\mathbf E}

\newcommand{\cd}{\mathcal D}

\newcommand{\cg}{\mathcal G}
\newcommand{\ch}{\mathcal H}

\newcommand{\cn}{\mathcal N}

\newcommand{\lp}{\left(}
\newcommand{\rp}{\right)}

\newcommand{\lcl}{\left\{}
\newcommand{\rcl}{\right\}}

\newtheorem{theorem}{Theorem}[section]

\newtheorem{corollary}[theorem]{Corollary}

\newtheorem{lemma}[theorem]{Lemma}

\newtheorem{proposition}[theorem]{Proposition}
\theoremstyle{remark}
\newtheorem{remark}[theorem]{Remark}

\def\tresp#1_#2^#3{\mathrel {\mathop{\kern 0pt#1}\limits_{#2}^{#3}}}

\begin{document}
\title{A model of continuous time polymer on the lattice}
\author[David M\'arquez, Carles Rovira, Samy Tindel]{David M\'arquez-Carreras \and Carles Rovira \and Samy Tindel}

\address{
{\it David M\'arquez-Carreras and Carles Rovira:}
{\rm Facultat de Matem\`{a}tiques,
Universitat de Barcelona, Gran Via 585, 08007 Barcelona, Spain}.
{\it Email: }{\tt davidmarquez@ub.edu, carles.rovira@ub.edu}.
\newline
$\mbox{ }$\hspace{0.1cm}
{\it Samy Tindel:}
{\rm Institut {\'E}lie Cartan Nancy, Nancy-Universit\'e, B.P. 239,
54506 Vand{\oe}uvre-l{\`e}s-Nancy Cedex, France}.
{\it Email: }{\tt tindel@iecn.u-nancy.fr}
}

\date{\today}

\begin{abstract}
In this article, we try to give a rather complete picture of the behavior of the free energy for a model of directed polymer in a random environment,
in which the polymer is a simple symmetric random
walk on the lattice $\Z^d$, and the environment is a collection $\{W(t,x);t\ge 0, x\in \Z^d\}$
of i.i.d. Brownian motions.
\end{abstract}

\subjclass[2000]{82D60, 60K37, 60G15}

\keywords{polymer measure, random environment, Gaussian estimates}

\maketitle

\section{Introduction}
\label{sec:intro} After two decades of efforts, the asymptotic
behavior of polymer measures, either in a discrete
\cite{AZ,Bo,CH,DSp,IS,Me,Peter,Pi} or continuous
\cite{BTV,CH2,CY,CY2,FH,RT} time setting, still remains quite
mysterious. Furthermore, referring to the articles mentioned
above, this problem is mainly tackled through the study of the
partition function of the measure. It is thus natural to try to
find a model for which a rather complete picture for the large
time behavior of the partition function is available. In the
current article, we shall show that one can achieve some sharp
results in this direction for a model of continuous time random
walk in a Brownian environment.

\smallskip

Indeed,
this paper is concerned with a model for a $d$-dimensional directed random walk
polymer in a Gaussian random medium which can be briefly
described as follows: the polymer itself, in the absence of any random
environment, will simply be modeled by a continuous time random walk
$\{b_t, t\ge 0\}$ on $\Z^d$. This process is
defined on a complete probability space $(\Omega_b, \mathcal{F},
\{P^x\}_{x\in \Z^d})$, where $P^x$ stands for the measure representing
the random walk starting almost surely from the initial
condition $x$; we write the corresponding expectation as $E_b^x$.
Let us recall that under $P^x$, the process is at  $x$ at time 0, it
stays there for an exponential holding time (with parameter
$\alpha=2 d$), and then jumps at one of the $2d$ neighbors
of $x$ in $\Z^d$ with equal probability. It stays there for an
exponential holding time independent of everything else and so on. Notice that a
given realization of $b$ belongs to the space
$\cd$ of paths $y:\R_+\to\Z^d$ such that
$y_t=\sum_{i\ge 0} x_i\1_{[\tau_i,\tau_{i+1})}(t)$, for a given sequence of increasing
positive times $(\tau_i)_{i\ge 0}$ and a given path $(x_i)_{i\ge 0}$ of a nearest neighbor
random walk $x$. We will denote by $\cd_{[0,t]}$ the restriction of such a space to
$\ot$.

\smallskip

The random environment in which our polymer lives is given by a family of
independent Brownian motions $\{ W(t,x);\ t \ge 0, x \in \Z^d\}$, defined
on some probability space $(\Omega, \mathcal{G}, {\bf P})$ independent
of $(\Omega_b, \mathcal{F},\{P^x\}_{x\in \Z^d})$. More specifically, this
means that $W$ is a centered Gaussian process satisfying
$${\bf E}\left[W(t,x) W(s,y)\right]=(s\wedge t) \delta(x-y),$$
where $\be$ denotes the expectation on $(\Omega, \mathcal{G}, {\bf
P})$, and $\delta$ stands for the discrete Dirac measure, i.e.
$\delta(z)=\1_{0}(z)$. We will also call $(\cg_t)_{t\ge 0}$ the
filtration generated by $W$.

\smallskip

Finally, our Gibbs type polymer measure is constructed as follows:
for any $t>0$, the energy of a given path  $b$ on
$[0,t]$ is given by the Hamiltonian:
$$
H_t(b):= \int_0^t W(ds,b_s).
$$
Observe that for a given path $b$, the last quantity has to be
understood as a Wiener integral with respect to $W$. In
particular, it is a centered Gaussian variable with variance $t$.
Now, for any $x \in \Z^d$, any $t\ge 0$ and a given $\beta>0$
which represents the inverse of the temperature, we define the
polymer measure   by the formula:
\begin{equation}\label{ezt}
d G^x_t(b)= \frac{e^{\beta H_t(b)}}{Z_t^x}\ dP^x(b),
\quad\mbox{ with }\quad
Z^x_t=E_b^x\left[e^{\beta H_t(b)}\right].
\end{equation}
The normalisation constant $Z^x_t$ above is called the partition function of the model.
In the sequel, we will also have to consider the
Gibbs average with respect to the polymer measure, defined as
follows: for any $t\ge 0$, $n\ge 1$, and for any bounded measurable
functional $f:\cd_{\ot}^n\to\R$, we set
$$
\langle f\rangle_t=\frac{E_b^x\left[f(b^1,\dots,b^n)
e^{\beta \sum_{l\le n} H_t(b^l)}\right]}{\left(Z_t^x\right)^n},
$$
where $\{b^1,\dots,b^n\}$ are understood as independent continuous time random
walks.

\smallskip

Let us say now a few words about the partition function: the first thing one can notice,
thanks to some invariance arguments, is that the asymptotic behavior of $Z_t^x$ does not
depend on the initial condition $x$. We shall thus consider generally $x=0$, and set $Z_t^0\equiv Z_t$. Moreover,
it can be checked (see Section \ref{sec:existence-free-energy} for further explanations)
that the following limits exist:
\begin{equation}\label{eq:def-p-beta}
p(\beta)=\lim_{t\uparrow\infty}
\frac{1}{t}\ \E \left[\log  Z_t \right]=
{\bf P}-{\rm a.s.}-\lim_{t\uparrow\infty} \frac{1}{t}\ \log  Z_t,
\end{equation}
and it can also be shown that $p(\beta)\le \beta^2/2$ by an elementary Jensen's type argument.
The quantity $p(\beta)$ is called the free energy of the system. It is then possible to
separate a region of {\it weak disorder} from a region of {\it strong disorder} according
to the value of $p(\beta)$, by saying that the polymer is in the strong disorder regime
if $p(\beta)< \beta^2/2$, and in the weak disorder regime otherwise. It should be mentioned
that this notion of strong disorder is rather called {\it very strong disorder} in \cite{CH2},
 the exact concept of strong disorder being defined thanks to some martingale
  considerations e.g. in \cite{CY,RT}. It is however believed that strong and
   very strong disorder coincide (see \cite{CH2} again). Furthermore, these
   notions have an interpretation in terms of
localization \cite{CH2} or diffusive behavior \cite{CY2} of the polymer.

\smallskip

With these preliminaries in hand, we will see now that some
sharp information on the partition function can be obtained for
the model under consideration. Namely, to begin with, the weak and
strong disorder regimes can be separated as follows:
\begin{proposition}\label{prop:weak-strong-disorder}
Let $Z_t\equiv Z_t^0$ be the normalization constant given by formula
(\ref{ezt}), and define a $\cg_t$-martingale $(M_t)_{t\ge 0}$ by $M_t=Z_t\exp(-\beta^2 t/2)$. Then:
\begin{enumerate}
\item
Whenever $d=1,2$ and $\beta>0$, we have $\lim_{t\to \infty}M_t=0$
in the $\bp$-almost sure sense, which means that the polymer is in the strong disorder regime.
\item
For $d\ge 3$ and $\beta$ small enough, the polymer is in the weak
disorder regime, i.e. $\lim_{t\to 0}M_t>0$, $\bp$-almost surely.
\item
For any dimension $d$ and $\beta>\beta_d$, the polymer is in the
very strong disorder regime, which means that $p(\beta)<\beta^2/2$.
\end{enumerate}
\end{proposition}
This kind of separation for the weak and strong disorder regime
has already been obtained for other relevant models, based on
discrete time random walks \cite{CH} or Brownian motions
\cite{CTV,CY2,RT}. However, the third point above can be sharpened
substantially, and the following almost exact limit holds true in
the continuous random walk context:
\begin{theorem}\label{thm:lim-free-energy}
Let $p(\beta)$ be the quantity defined at (\ref{eq:def-p-beta}), and
$\varepsilon_0$ be  a  given arbitrary positive constant. Then, there
exists $\beta_0=\beta_0(d)>0$ such that
$$C_0 \frac{\beta^2}{\log\left(\beta\right)}\ (1-\varepsilon_0)
\le p(\beta) \le C_0 \frac{\beta^2}{\log\left(\beta\right)}\
(1+\varepsilon_0), \qquad {\rm for}\ \beta \ge \beta_0,$$
where $C_0$ is a strictly positive
positive constant which will be defined by relation (\ref{eq:def-C0}).
\end{theorem}
Putting together Proposition \ref{prop:weak-strong-disorder} and
Theorem \ref{thm:lim-free-energy}, we thus get a remarkably precise
picture as far as the free energy of the system is concerned.
\begin{remark}
Many of our results would go through without much effort for a wide
class of spatial covariance of the medium $W$, as done in \cite{CTV}.
We have sticked to the space-time white noise case in the current article for sake of simplicity.
\end{remark}

\smallskip

Our paper is divided as follows: at Section \ref{sec:existence-free-energy},
we recall some basic facts about the partition function of the polymer model.
Section \ref{sec:limit-free-energy} is the bulk of our article, and is devoted
to a sharp study of the free energy in the low temperature region, along the
lines of the Lyapunov type result \cite{CKM,CV,KVV}.
At Section \ref{sec:weak-strong-disorder}, the first two items of
Proposition \ref{prop:weak-strong-disorder} are shortly discussed.

\section{Basic properties of the free energy}
\label{sec:existence-free-energy}

Since it will be essential in order to show Theorem \ref{thm:lim-free-energy},
we will first devote the current section to show briefly that $Z_t$ converges
almost surely to a constant $p(\beta)$, which can be done along the same lines
as in \cite{RT}. First of all, some standard
arguments yield the following asymptotic result:
\begin{proposition} \label{pmarkov}
For $t>0$, define
$
p_t(\beta)= \frac{1}{t}\ \E \left[\log Z_t^x\right].
$
Then, for all $\beta>0$, there exists a constant
$p(\beta)>0$ such that
$$p(\beta)\equiv \lim_{t\uparrow\infty} p_t(\beta)=\sup_{t\ge 0}
p_t(\beta).
$$
Furthermore, $p(\beta)\le \beta^2/2$.
\end{proposition}

\begin{proof}
It can be proved e.g. as in \cite{RT}. More specifically, we should first show a Markov
decomposition for $Z_t^x$ as in \cite[Lemma 2.4]{RT}; then we can argue as in
\cite[Proposition 2.5]{RT} in order to get the announced limit for $p_t(\beta)$.
The bound $p(\beta)\le \beta^2/2$ can be checked
using Jensen's inequality.

\end{proof}

\begin{remark}
Due to the spatial homogeneity of $W$, the above limit does not depend
on $x\in\Z^d$. Hence, from now on,  we will choose $x=0$ for our computations.
Furthermore,
in the sequel, when $x=0$ we will write $Z_t$ and $E_b$ instead of
$Z_t^0$ and $E_b^0$, respectively.
\end{remark}

In order to get the almost sure convergence of $p_t(\beta)$, we will need some
concentration inequalities which can be obtained by means of Malliavin calculus tools.
Let us briefly recall here the main features of this theory, borrowed from \cite{Nu}.
First of all, let us notice that our
Hamiltonian $H_t(b)$ can be written as follows:
$$H_t(b):= \int_0^t W(ds,b_s)=\int_0^t \sum_{x\in \Z^d} \delta_x(b_s) W(ds,x).$$
This leads to the following natural definition of an underlying Wiener space in our
context:
set $\mathcal{H}=L^2(\mathbb{R}_+ \times \Z^d)$, endowed with the norm
$|h|_\mathcal{H}^2=\int_{\mathbb{R}_+}  \sum_{x\in \Z^d} |h(t,x)|^2 dt.$
Then there exists a zero-mean isonormal Gaussian family
$\{W(h); h\in
\mathcal{H}\}$ defined by
$$
W(h)=\int_{\mathbb{R}_+} \sum_{x\in\Z^d} h(t,x)  W(dt,x),
$$
and it can be shown that $(\Omega,\ch,\bp)$ is an abstract Wiener space.
Denote now by $\mathcal{S}$ the set of
smooth functionals defined on this Wiener space, of the form
$$
F=f(W(h_1),\dots,W(h_k)),
\quad \mbox{ for }\quad
k\ge 1,\ h_i\in \mathcal{H},\ f\in C_b^\infty(\R^k).
$$
Then its Malliavin derivative is defined as
$$D_{t,x}F=\sum_{i=1}^k \partial_i f(W(h_1),\dots, W(h_k))
h_i(t,x).$$ As usual, the operator $D: \mathcal{S}\longrightarrow
\mathcal{H}$ is closable and we can build the family of Sobolev spaces
$\mathbb{D}^{1,p}$, $p\ge 1$, obtained by completing $\mathcal{S}$
with respect to the norm
$\|F\|_{1,p}^p= \E\left[|F|^p\right]+\E
\left[|DF|_\mathcal{H}^p\right].$   The following chain
rule is also available for $F\in \mathbb{D}^{1,p}$: if $\psi: \mathbb{R}
\longrightarrow \mathbb{R}$ is a smooth function, then $\psi(F)
\in \mathbb{D}^{1,q}$ for any $q<p$ and
\begin{equation}\label{eq:chain-rule}
D \psi(F)=\psi'(F) DF.
\end{equation}

\smallskip

We are now ready to prove the almost sure limit of $\log(Z_t)/t$.
\begin{proposition}
We have that
$$\bp- {\rm  a.s.}-\lim_{t\uparrow\infty}
\frac{1}{t}\ \log  Z_t=p(\beta).$$
\end{proposition}

\begin{proof}
It is easily shown that $Z_t\in\D^{1,2}$, and
by differentiating in the Malliavin calculus sense, we obtain, if $s<t$:
$$D_{s,x} Z_t= \beta E_b\left[ ( D_{s,x} H_t(b)) e^{\beta
H_t(b)}\right]=\beta E_b\left[ \delta_x(b_s)  e^{\beta
H_t(b)}\right].$$
Thus, if $U_t=\frac{1}{t}\ \log Z_t$, we have, thanks to the chain rule (\ref{eq:chain-rule}),
$$D_{s,x} U_t=\frac{D_{s,x} Z_t}{t Z_t}=\left\{\begin{array}{ll}
 \beta \ \frac{\langle
\delta_x(b_s)\rangle_t}{t},\ & s\leq t,\\[2mm]
0,& s>t.\end{array}\right.$$
So,
$$|D U_t|^2_\mathcal{H}= \frac{\beta^2}{t^2}\
\int_0^t \sum_{x\in \Z^d} \langle \delta_x(b_s^1)
\delta_x(b^2_s)\rangle_t ds =\frac{\beta^2}{t^2}\ \int_0^t \langle
\delta_0(b_s^1-b^2_s)\rangle_t ds\le \frac{\beta^2}{t} .$$ Now,
since $|D U_t|^2_\mathcal{H}$ is bounded and tends to 0 as $t\to\infty$, we
can prove the almost sure limit by means of a concentration
inequality (see, for instance, \cite[Proposition 2.1]{RT}) and a
Borel-Cantelli type argument.

\end{proof}

\section{Exact limit for the free energy at low temperature}
\label{sec:limit-free-energy}

The aim of this section is to show our Theorem \ref{thm:lim-free-energy},
by means of some Gaussian tools which have been already used for various
models of polymers \cite{CTV} or stochastic PDEs \cite{CKM,CM,FV}. It should
be mentioned at this point that the reference \cite{CKM} is especially relevant
for us: in fact, our aim here is to clarify some of the arguments therein, and
adapt them to the polymer context at the same time.

\subsection{Strategy}
\label{sec:strategy}
In order to understand how the free energy will be computed, let us introduce first
some additional useful notations:
let $P_n$
be the set of paths of a discrete time random walk of length $n$
starting at $0$, and $S_{n,t}$ be the set of possible
times of the jumps of the continuous time random walk $b$ in $\ot$, namely:
$$S_{n,t} = \{\bs =(s_0,\dots,s_n);\ 0=s_0 \le s_1 \le
\ldots \le s_n \le t \}.$$
For $\bs  \in S_{n,t}$ and
$\bx=(x_0,\dots,x_n) \in P_n$, we also define
$$B_n(\bs ,\bx)=
\sum_{j=0}^n \left[W(s_{j+1},x_j)- W(s_{j},x_j)\right].$$
For any positive $t$, let $N_t$ be the number of jumps of $b$ in $\ot$, which
is known to be a Poisson process with intensity $2d$. Then one can decompose the
Hamiltonian $H_t(b)$ according to the values of $N_t$ in order to obtain:
\begin{eqnarray}
Z_t &=& \sum_{n=0}^{\infty} E_b\left[ e^{\beta\int_0^t
 W(ds,b_{s})}\ \1_{\{N_t=n\}} \right]\nonumber\\
&=& \sum_{n=0}^{\infty}\frac{(2 d  t)^n}{n!}\ e^{-2 d  t} \
E_b\left[ e^{\beta \int_0^t
 W(ds,b_{s})} | N_t=n \right]\nonumber\\
&=& \sum_{n=0}^{\infty} e^{-2 d  t}\sum_{\bx\in P_n}
\int_{S_{n,t}} e^{\beta B_n(\bs ,\bx)
 }  ds_1\cdots ds_n.\label{e5}
\end{eqnarray}

\smallskip

With these preliminaries in mind, we can now sketch the strategy we shall follow
in order to obtain the lower and upper bounds  on $p_t(\beta)$ announced at Theorem
\ref{thm:lim-free-energy}. Indeed, the basic idea is that one should find an equilibrium between two constraints:

\smallskip

\noindent
{\it (i)} The more the random walk jumps, the more it will be able to see the peaks of the
energy, represented by $$\sup_{j\le n,\bs\in S_{j,t},\bx\in P_j}B_j(\bs ,\bx).$$ We shall
see that, roughly speaking, $\sup_{j\le r t,\bs\in S_{j,t},\bx\in P_j}B_j(\bs ,\bx)$ is of
order $r^{1/2}t$ for a given $r>0$.

\smallskip

\noindent
{\it (ii)} The number of jumps of the random walk has an entropy cost, which
is represented in formula (\ref{e5}) by the area of $S_{n,t}$. It is a well
known fact that this area decreases as $n!$.

\smallskip

\noindent After reducing our calculations to this optimization
problem, it will be easily seen that the accurate choice for $r$
is of order $(\beta/\log(\beta))^2$, a fact which has already been
outlined in \cite{CTV}. This means that, under the influence of
the environment when $\beta$ is large enough, the random walk is
allowed to jump substantially more -the typical number of jumps
before $t$ is of order $(\beta/\log(\beta))^2 t$- than in the free
case (for which this typical number is of order $2 d t$). By
elaborating this kind of considerations, we shall obtain our bound
$C_0 \beta^2/ \log(\beta)$ for $p(\beta)$.

\smallskip

We are now ready to perform our first technical step, which is the control on the Gaussian random field $B_n$.

\subsection{Control of the random field $B_n$}
\label{sec:control-B-n}
Let us mention that,
in order to control the supremum of the random Gaussian fields
we will meet, we shall use the following classical results, borrowed from \cite{A}.
\begin{theorem}\label{the21}
Let $G(t)$ be a mean zero Gaussian field over a set $T$ with the
associated pseudo-metric
$$d(t_1,t_2)=\left(\E\left[G(t_1)-G(t_2)\right]^2\right)^{1/2},$$
and assume that the metric space $(T,d)$ is compact. Let $\cn (\eta)$
be the metric entropy associated with
$d$, i.e. the smallest number of closed $d$-balls of radius
$\eta$ needed to cover $T$. Then, there exists a universal
constant $K$ such that
$$\E \sup_T G(t)\le K \int_0^\infty \sqrt{\log \cn (\eta)}\
d\eta,$$ provided that the right hand side is finite.
\end{theorem}
This useful theorem for the computation of the mean of $\sup_{t\in T} G(t)$ has
generally to be completed by a control on the fluctuations of Gaussian fields, given by the following result:
\begin{theorem}\label{the22}
Let $G(t)$ be a mean zero Gaussian field over a set $T$, and suppose that the
sample paths of $G$ are bounded almost surely. Then, we have
$\E\sup_T G(t)< \infty,$
and for all \ $y>0$,
$${\bf P} \left\{\left| \sup_T G(t)- \E\sup_T G(t)\right|>y\right\}\le
2 e^{-y^2/2\tau^2}\ ,$$ where $\tau^2=\sup_T \E [G(t)^2]$.
\end{theorem}

\medskip

\noindent
With these Gaussian tools in hand,
we can now  control $B_n$ as follows:
\begin{proposition}\label{pentro}
For $r\ge 0$, let $T_r$ be the space of paths of the
continuous time random walk starting at $x=0$ with
no more than $rt$ jumps. Then, we have
\begin{equation}\label{eboubd}\E\sup_{(n,\bs ,\bx)\in T_r}
B_n(\bs ,\bx) \le C \sqrt{[r t]
t},\end{equation} for some positive constant $C$, where $[u]$ stands
for the integer part of a real number $u$. Moreover,
\begin{equation}\label{emilu}
\lim_{t\rightarrow +\infty} \frac{1}{t}\
\E\sup_{(n,\bs ,\bx)\in T_r}
B_n(\bs ,\bx)\equiv F(r) <+\infty.\end{equation}
Finally, the following scaling identity holds true for the function $F$: for any $r>0$,
we have $F(r)=\sqrt{r}\ F(1)$.
\end{proposition}
\begin{remark}
Observe that, in order to describe an element of $T_r$, we have to know the number
of jumps $n$, the times of jumps $\bs \in S_{n,t}$ and
the paths $\bx\in P_n$. Hence, the family
$B_n(\bs ,\bx)$ can be considered as a Gaussian field over
$T_r$. Moreover, the random variable
$\sup_{(n,\bs ,\bx)\in T_r}B_n(\bs ,\bx)$ and its expectation only depend
on the parameters $r$ and $t$.
\end{remark}

\medskip

\noindent  Before giving the proof Proposition \ref{pentro}, let us also state the
following corollary, which can be proved along the same lines.
\begin{corollary}\label{centro}
Let $T_r^\varrho$ be the space of paths of the continuous time
random walk starting at $x=0$ with no more than
$rt$ jumps,
with the additional constraint that the jumps are separated from each other and from the
endpoints of the interval $[0,t]$ by a distance of at least
$2\varrho$. Then, the limit
\begin{equation*}\lim_{t\rightarrow +\infty} \frac{1}{t}\
\E\sup_{(n,\bs ,\bx)\in T^\varrho_r} B_n(\bs ,\bx)\equiv
F^\varrho(r)<+\infty
\end{equation*}
exists and, moreover,  $F^\varrho(r)=\sqrt{r}\ F^{r\varrho}(1)$.
\end{corollary}

\begin{proof}[Proof of Proposition \ref{pentro}:]
 In order to
obtain the bound (\ref{eboubd}) we will use Theorem \ref{the21},
which involves the entropy of the Gaussian field
$B_n(\bs ,\bx)$ over $T_r$. Let us then estimate this entropy: for $n\ge 0$,
$\bs , \bsp\in S_{n,t}$ and
$\bx\in P_n$, we define the distance between
$(n,\bs ,\bx)$ and
$(n,\bsp,\bx)$ as
$$d((n,\bs ,\bx),(n,\bsp,\bx))
=\sqrt{\E\left[|B_n(\bs ,\bx)-B_n(\bsp,\bx)|^2\right]};$$
and let $\cn (\eta)$ be the entropy function related with
$B_n(\bs ,\bx)$ and $T_r$. Since
$\E\left[B_n(\bs ,\bx)^2\right]=t$, the
diameter of $T_r$ is smaller than $2 \sqrt{t}$. Assume then that
$\eta\le  2 \sqrt{t}$. It is not difficult to check (see \cite{FV} for a similar computation) that
$$d^2((n,\bs ,\bx),(n,\bsp,\bx))\le
2\sum_{j=0}^n \left|s_j -s'_j\right|.$$
Thanks to this identity, one can construct a $\eta$-net in $T_r$ for the pseudo-metric $d$.
 It is simply based on all the paths of the random walk with any position vector $\bx$, and j
 ump times in the subset
$\hat S_{n,t}$ of $S_{n,t}$ defined as  follows: $\hat S_{n,t}$ is the set of elements
$\bs =(s_0,\dots,s_n)\in S_{n,t}$ where all $s_j$ are
integer multiples of $\eta^2 (4n)^{-1}$ (notice that some of the $s_j$'s can be equal). It
is readily checked that all these paths form a  $\eta$-net in $T_r$, and furthermore,
the cardinal of $\hat S_{n,t}$ can be bounded as:
$$\sharp \hat S_{n,t}=
\frac{1}{n!}\ \left[\frac{4nt}{\eta^2}\right]^n
\le \left(\frac{C
t}{\eta^2}\right)^n,
$$
where in the last step we have used the inequality $n!\ge
(n/3)^n$. So, owing to the fact that $\sharp P_n$ is bounded by
$(2d)^n$, that $n\le [rt]$ and $\eta\le  2 \sqrt{t}$, we obtain
that
\begin{equation}\label{eentropy}
\cn (\eta)\le \sum_{n=0}^{[rt]} (2d)^n \left(\frac{C
t}{\eta^2}\right)^n
 \le \left(\frac{C t}{\eta^2}\right)^{[rt]}
\le 1 \vee \left(\frac{C
t}{\eta^2}\right)^{[rt]}.
\end{equation}
Thus, using Theorem
\ref{the21}, we end up with:
\begin{multline*}
\E\sup_{(n,\bs ,\bx)\in T_r} B_n(\bs ,\bx) \le C \int_0^{+\infty}
\sqrt{\log\cn (\eta)}\ d\eta \le C \sqrt{[rt]}\int_0^{2\sqrt{t}}
\sqrt{\log\left(\frac{C t}{\eta^2}\right)}\
d\eta  \\
\le  C \sqrt{[rt] t}\int_0^{1} \sqrt{-\log \eta^2}\ d\eta
\le C \sqrt{[r t] t}.
\end{multline*}
Our claim (\ref{emilu}) is now easily
verified thanks to a super-additive argument and, since $W$ is a
Wiener process in $t$, we also get, by Brownian scaling:
$$
F(r)=\lim_{t\rightarrow +\infty} \frac{1}{t}\
\E\sup_{(n,\bs ,\bx)\in T_r}
B_n(\bs ,\bx)= \sqrt{r
}\lim_{t\rightarrow +\infty} \frac{1}{t}\
\E\sup_{(n,\bs ,\bx)\in T_1}
B_n(\bs ,\bx)=\sqrt{r}
F(1).
$$

\end{proof}

\smallskip

\noindent
Notice that Corollary \ref{centro} depends on the discretization type parameter
$\varrho$, which is useful for technical purposes. However,
in order to get our final estimate on $p(\beta)$, we will need the
following lemma, which relates $F^\varrho(1)$ and $F(1)$.
\begin{lemma}\label{lentro} With the notations of Proposition \ref{pentro} and
Corollary \ref{centro} we have that
$$F^\varrho(1)\quad  \tresp
\longrightarrow_{\varrho\rightarrow 0}^{} \quad F(1).$$
\end{lemma}
\begin{proof}
Let us denote by $\phi$ and $\phi^{\varrho}$ the quantities defined by:
\begin{equation*}
\phi^\varrho(r,t)=\E\sup_{(n,\bs ,\bx)\in T_r^\varrho} B_n(\bs ,\bx),
\quad\mbox{and}\quad
\phi(r,t)=\E\sup_{(n,\bs ,\bx)\in T_r} B_n(\bs ,\bx).
\end{equation*}
We will prove that, for an arbitrary
$\varepsilon>0$, we can find $\varrho_0>0$ such that
\begin{equation}\label{efity}
F(1)- \lim_{t\rightarrow \infty}
\frac{\phi^\varrho(r,t)}{t}<\varepsilon,\qquad {\rm for}\
\varrho\le \varrho_0.
\end{equation}
Indeed, the triangular inequality implies that
\begin{equation}\label{epero}
\left|F(1)-\frac{\phi^\varrho(r,t)}{t}\right|\le
\left|F(1)-\frac{\phi(r,t)}{t}\right|+\left|\frac{\phi(r,t)}{t}-\frac{\phi^\varrho(r,t)}{t}\right|.\end{equation}
Now, according to Proposition \ref{pentro}, there exists $t_0$ such that, for any
$t\ge t_0$,
\begin{equation}\label{epero2}
\left|F(1)-\frac{1}{t}\ \phi(r,t)\right|\le
\frac{\varepsilon}{2}.\end{equation}
Furthermore, for a fixed $t_0$, since
$\phi^\varrho(r,t_0)\rightarrow \phi(r,t_0)$ when $\varrho
\rightarrow 0$, there exists $\varrho_0>0$ such that
\begin{equation}\label{epero3}
\left|\frac{1}{t_0}\left[\phi(r,t_0)-\phi^\varrho(r,t_0)\right]\right|\le
\frac{\varepsilon}{2}\qquad {\rm for}\ \varrho\le
\varrho_0.
\end{equation}
Finally, a super-additivity type argument easily yields the fact that $\phi^\varrho(r,t)/t$ is
increasing in $t$. This property, together with (\ref{epero})-(\ref{epero3}) implies
(\ref{efity}), which ends our proof.

\end{proof}

\smallskip

Let us now complete the information we have obtained on the expected value
of $\sup B_n(\bs ,\bx)$ by a study of the fluctuations in $\bs$ of the field
$B_n(\bs ,\bx)$. To this purpose, let us introduce a little more notation:
for $r,\rho>0$, let $Y_{r,\rho}\subseteq T_r\times T_r$ be the set defined as:
\begin{multline}\label{eq:def-y-r-rho}
Y_{r,\rho}=
\Big\{
((n,\bs ,\bx),(n',\bsp,\bxp ))\in T_r\times T_r;\,
n=n'\le [rt], \, \bx=\bxp, \\
|s_j-s'_j|\le \rho \mbox{ for } 1\le j \le n
\Big\}.
\end{multline}

\begin{proposition}\label{plult}
For  $r,\rho>0$, let $\Upsilon(r,\rho)=\limsup_{t\rightarrow
\infty} \frac{1}{t}\ \E \sup_{Y_{r,\rho}} A_n(\bs ,\bsp,\bx),$
where $A_\cdot(\cdot,\cdot,\cdot)$ is  a fluctuation Gaussian
field defined on $Y_{r,\rho}$ as follows:
\begin{equation}\label{edefan}
A_n(\bs ,\bsp,\bx)=
B_n(\bs ,\bx)-B_n(\bsp,\bx).\end{equation}
Then, $\Upsilon(r,\rho)=\sqrt{r}\ \Upsilon(1,r\rho)$
and
\begin{equation}\label{elimupsil}
\lim_{\rho\rightarrow 0}\Upsilon(r,\rho)=0.
\end{equation}
\end{proposition}
\begin{proof}
The scaling property $\Upsilon(r,\rho)=\sqrt{r}\ \Upsilon(1,r\rho)$ is
easily shown. Let us concentrate then on relation (\ref{elimupsil}): The bound
(\ref{eentropy}) implies that the entropy for the field
$B_n(\bs ,\bx)$ over $T_1$ is bounded as
follows: \begin{equation}\label{ebon12}\cn _B(\eta)\le 1
\vee \left(\frac{C t}{\eta^2}\right)^{t}.\end{equation} Hence,
the entropy for the field $A_n(\bs ,\bsp,\bx)$ over
$Y_{1,\rho}$ also satisfies:
\begin{equation}\label{ebon23}\cn _A(\eta)\le
\cn ^2_B\big(\frac{\eta}{2}\big).\end{equation} The second
ingredient to start the proof of (\ref{elimupsil}) is a bound on the
diameter of $Y_{1,\rho}$ in the canonical metric associated to $B$:
\begin{equation}\label{ebon22}
{\rm Diam}(Y_{1,\rho})\le
2 \sup_{Y_{1,\rho}}  \E^{1/2}\left[A_n^2(\bs ,\bsp,\bx)\right]
\le 2 \sqrt{\rho t}.\end{equation} Then, Theorem \ref{the21}
together with (\ref{ebon12}), (\ref{ebon23}) and (\ref{ebon22})
yield
\begin{eqnarray*}
\E \sup_{Y_{r,\rho}} A_n(\bs ,\bsp,\bx)
&\le&  C
\int_0^{2\sqrt{\rho t}} \sqrt{\log \cn _A(\eta)} \ d\eta
\le  C \int_0^{2\sqrt{\rho t}} \sqrt{\log
\left(\frac{4Ct}{\eta^2}\right)^{2t}} \ d\eta \\
&=& 2C \sqrt{2t}\int_0^{2\sqrt{\rho t}} \sqrt{\log
\left(\frac{4Ct}{\eta^2}\right)} \ d\eta
\le  Ct
\int_0^{\sqrt{\rho/C}} \sqrt{-2\log \zeta} \ d\zeta.
\end{eqnarray*}
The proof is now finished along the same lines as for Proposition \ref{pentro}.

\end{proof}

\smallskip

It is worth mentioning at this point that the function $F$ emerging
at relation (\ref{emilu}) is the one which allows us to define the
constant $C_0$ in Theorem \ref{thm:lim-free-energy}. Indeed, this constant is simply given by
\begin{equation}\label{eq:def-C0}
C_0=  \frac18 F(1)^2.
\end{equation}
With these preliminaries in hand, we are now ready to proceed to the
proof of our Theorem \ref{thm:lim-free-energy}. This proof will be divided
between the lower and the upper bound, for which
the techniques involved are slightly different.

\subsection{Proof of Theorem \ref{thm:lim-free-energy}: the lower bound}

Recall that we wish to prove that, for $\beta$ large enough and an arbitrary
positive constant $\varepsilon_0$, we have
$p(\beta)\ge C_0 \frac{\beta^2}{\log\left(\beta\right)}\ (1-\varepsilon_0)$.
Now, since $p(\beta)$ exists
and is non-random, we  only need to prove that, for $\beta \ge
\beta_0$ and for $t$ large enough,
\begin{equation}\label{elower1}
{\bf P} \left\{\frac{1}{t \beta^2}\ \log
\left[\sum_{n=0}^{[r(\beta)t]} e^{-2dt}\sum_{\bx \in
P_n} \int_{S_{n,t}} e^{ \beta
B_n(\bs ,\bx )
 }  ds_1\cdots ds_n\right]\ge \frac{C_0(1-\varepsilon_0)}{\log
 (\beta)}\right\}\ge \frac12,\end{equation}
where, as mentioned at Section \ref{sec:strategy}, the parameter $r(\beta)$ will be chosen as:
$$r(\beta)=\frac{C_0}{2} \frac{\beta^2}{\log^2(\beta)}.$$

\smallskip

\noindent
{\it Step 1: Reduction of the problem.}
Observe that, on the one hand, we have
 \begin{equation}\label{elower2}\begin{array}{l}
\displaystyle \frac{1}{t \beta^2}\ \log
\left[\sum_{n=0}^{[r(\beta)t]} e^{-2 d t}\sum_{\bx \in P_n} \int_{S_{n,t}} e^{\beta
B_n(\bs ,\bx )
 }  ds_1\cdots ds_n\right]\\[4mm]
 \qquad \displaystyle \ge \frac{1}{t \beta^2}\ \log
\left[\sup_{n\le [r(\beta)t]} e^{-2 d t}\sum_{\bx \in
P_n} \int_{S_{n,t}} e^{\beta B_n(\bs ,\bx )
 }  ds_1\cdots ds_n\right]\\[4mm]
 \qquad \displaystyle \ge - \frac{2
d}{\beta^2} + \frac{1}{t\beta^2}\ \log
\left[\sup_{n\le [r(\beta)t]} \sup_{\bx \in P_n}
\int_{S_{n,t}} e^{\beta B_n(\bs ,\bx )
 }  ds_1\cdots ds_n\right].\end{array}\end{equation}
On the other hand, for $\beta\ge\beta_1$ and  $\beta_1$ large enough, we have
$\frac{2d}{\beta^2} \le C_0 \frac{\varepsilon_0}{2} \frac{1}{ \log(\beta)}$. Furthermore, since $n \le \frac{C_0}{2} \frac{ \beta^2}{ \log^2(\beta)} t$, we obtain:
$$
\frac{ 2n \log(\beta)}{ t \beta^2} \le  C_0 \frac{1}{\log(\beta)},
\quad\mbox{and hence}\quad
\frac{ \log(\beta^{-2n})}{ t \beta^2} \ge -  C_0
\frac{1}{\log(\beta)}.
$$
Plugging this inequality into (\ref{elower2}), we get that in order to show
(\ref{elower1}), it is sufficient to prove that, for $\beta,t$ large enough,
the following relation holds:
\begin{equation}\label{elower3}
{\bf P} \left\{\frac{1}{t \beta^2}\ \log \left[\sup_{n\le
[r(\beta)t]} \sup_{\bx \in P_n} \int_{S_{n,t}}
\beta^{2n} e^{\beta B_n(\bs ,\bx )
 }  ds_1\cdots ds_n\right]\ge C_0(2-\varepsilon_0/2) \frac{1}{\log
 (\beta)}\right\} \ge \frac12.\end{equation}

\smallskip

In order to show relation (\ref{elower3})
set first  $\hat r(\beta) = (2 \beta^2)^{-1}$.
For $r,\varrho>0$, define also a set $\hat Y_{r,\varrho}$ by
\begin{equation}\label{eq:def-hat-Y}
\hat Y_{r,\varrho}=
\lcl
((n,\bs ,\bx ),(n,\bsp ,\bx ))\in Y_{r,\varrho}; \,
|s_{j+1}-s_j|\ge 2\varrho, \mbox{ for } j=1\ldots n
\rcl,
\end{equation}
where $Y_{r,\varrho}$ is defined by (\ref{eq:def-y-r-rho}).
Finally, for
$((n,\bs ,\bx ),(n,\bsp ,\bx ))\in
Y_{r(\beta),\hat r(\beta)}$, we set
$$\xi_{\beta} (\bsp )=\{\bs ;\ |s_j-s_j'|\le
\hat r(\beta),\ 1\le j\le n\}.
$$
Using the definition of these sets, and recalling that the field
$A_\cdot(\cdot,\cdot,\cdot)$ has been defined by (\ref{edefan}), we
have that, for any $(n,\bsp ,\bx )\in T_{r(\beta)}^{2\hat r(\beta)}$,
\begin{eqnarray}
\int_{S_{n,t}} \beta^{2n} e^{\beta
B_n(\bs ,\bx )
 }  ds_1\cdots ds_n&\ge & \beta^{2n} e^{\beta B_n(\bsp ,\bx )
 }\int_{S_{n,t}} e^{\beta A_n(\bs , \bsp , \bx )
 }  ds_1\cdots ds_n\nonumber\\
 &\ge & \beta^{2n} e^{\beta B_n(\bsp ,\bx )
 }\int_{\xi_{\beta}(\bsp )} e^{\beta A_n(\bs , \bsp , \bx )
 }  ds_1\cdots ds_n\label{eblip}\\
&\ge &e^{\beta B_n(\bsp ,\bx )
 }\ \exp\left(-\sup\{ \beta A_n(\bs , \bsp ,
 \bx );\ \bs \in
 \xi_{\beta} (\bsp )\}\right).\nonumber
 \end{eqnarray}
Since (\ref{eblip}) is true for any
$(n,\bsp ,\bx )\in T_{r(\beta)}^{2\hat r(\beta)}$,
and with (\ref{elower3}) in mind, we obtain:
\begin{multline}\label{etamber}
\frac{1}{t \beta^2}\ \log \left[\sup_{n\le
[r(\beta)t]} \sup_{\bx \in P_n} \int_{S_{n,t}}
\beta^{2n} e^{\beta B_n(\bs ,\bx )
 }  ds_1\cdots ds_n\right]  \\
 \ge \frac{1}{t \beta}\
\sup_{T_{r(\beta)}^{2\hat r(\beta)}}B_n(\bsp ,\bx )
-\frac{1}{t \beta }\ \sup_{\hat Y_{r(\beta),\hat
r(\beta)}}A_n(\bs , \bsp ,
 \bx ).
\end{multline}

\smallskip

\noindent
{\it Step 2: Study of the term involving} $B_n.$
Let us consider $\varepsilon_1 >0$.
According to Corollary \ref{centro}, there exists a constant $\beta_2>0$ such that
for $\beta\ge \beta_2$, we can find $t'=t'(\varepsilon_1,\beta)$ satisfying, for $t\ge t'$:
\begin{eqnarray*}
\frac{1}{\beta} \E \left[ \sup_{T_{r(\beta)}^{2\hat r(\beta)}}
B_n(\bsp ,\bx ) \right] &\ge & \frac{t}{\beta} \left[F^{2\hat
r(\beta)}(r(\beta))-\frac{\varepsilon_1
\sqrt{ r(\beta)}}{2}\right]\\
&=& \frac{t}{\beta} \left[ \sqrt{r(\beta)}F^{2\hat
r(\beta)r(\beta)}(1)-\frac{\varepsilon_1
\sqrt{r(\beta)}}{2}\right]\\
&= & \frac{t}{\beta} \sqrt{r(\beta)} \left[F^{2\hat
r(\beta)r(\beta)}(1)-\frac{\varepsilon_1}{2}
\right].\end{eqnarray*} Notice that we have chosen $r(\beta)$ and
$\hat r(\beta)$ so that
\begin{equation}\label{limit0}
\lim_{\beta \to \infty} r(\beta) \hat r(\beta) =0.
\end{equation}
Hence, applying now Lemma \ref{lentro}, we get that there exists $\beta_3>0$ such that
for any $\beta \ge \beta_3$, it holds that
$$
\frac{1}{\beta} \E \left[ \sup_{T_{r(\beta)}^{2\hat r(\beta)}}
B_n(\bsp ,\bx ) \right] \ge \frac{t}{\beta} \sqrt{r(\beta)}
\left[F(1)-\varepsilon_1 \right].
$$
Thus,
 Theorem \ref{the22} for
$\tau^2= \frac{t}{\beta^2}$ and $y=\frac{t}{\beta}\varepsilon_1
\beta \sqrt{r(\beta)}$ implies, for $t$ large enough,
\begin{equation}\label{esi16}
{\bf P}\left\{ \sup_{T_{r(\beta)}^{2\hat r(\beta)}} \frac{1}{\beta}
B_n(\bsp ,\bx ) < \frac{t}{\beta}
\sqrt{r(\beta)} \left[F(1)-2\varepsilon_1 \right]\right\}\le 2
\exp\left\{-\frac{t \varepsilon_1^2 {r(\beta)}}{2}\right\}\le
\frac14.\end{equation}

\smallskip

\noindent
{\it Step 3: Study of the term involving} $A_n.$
Invoking Proposition \ref{plult}, we have
\begin{equation*}\label{ekapp3}
\limsup_{t \rightarrow \infty} \frac{1}{t}\ \E
\sup_{Y_{r(\beta),\hat r(\beta)}}A_n(\bs ,
\bsp ,
 \bx )< \Upsilon (r(\beta),\hat r(\beta)) + \frac{\varepsilon_1}{2}.\end{equation*}
 Then, since (\ref{limit0}) holds true, and using (\ref{elimupsil}), we can
choose $\beta_4>0$ such that for $t$ large enough and $\beta \ge
\beta_4$ we have:
 $$ \frac{1}{\beta} \E \left[ \sup_{Y_{r(\beta),\hat r(\beta)}} \beta A_n(\bs , \bsp ,
 \bx ) \right]\le   \frac{t }{\beta} \sqrt{r(\beta)}
 \left( \Upsilon(1,r(\beta)\hat r(\beta)) + \frac{\varepsilon_1}{2} \right)
 \le \varepsilon_1 t \beta \sqrt{r(\beta)}.$$
Hence, Theorem \ref{the22} for $\tau^2\le 2 \frac{t}{\beta^2}$ and
$y=\varepsilon_1 \frac{t }{\beta} \sqrt{r(\beta)}$ yields, for $t$
large enough,
\begin{equation}\label{esi18}
{\bf P}\left\{ \sup_{Y_{r(\beta),\hat r(\beta)}} \frac{1}{\beta}
A_n(\bs , \bsp ,
 \bx ) > 2 \varepsilon_1 \frac{t }{\beta} \sqrt{r(\beta)}
\right\}\le 2 \exp\left\{-\frac{t \varepsilon_1^2
r(\beta)}{4}\right\}\le \frac14.\end{equation}

\smallskip

\noindent
{\it Step 4: Conclusion.}
Choose $\beta_0=\beta_1 \vee \beta_2 \vee \beta_3  \vee \beta_4$,
and consider $\beta\ge \beta_0$. Define also the sets:
\begin{eqnarray*}
\Omega_1&=&\left\{ \sup_{Y_{r(\beta),\hat r(\beta)}}
\frac{1}{\beta} A_n(\bs , \bsp ,
 \bx ) > 2 \varepsilon_1 \frac{t}{ \beta} \sqrt{r(\beta)}
\right\},\\
\Omega_2&=&\left\{ \sup_{T_{r(\beta)}^{2\hat r(\beta)}} \frac{1}{\beta}
B_n(\bsp ,\bx ) < \frac{t}{ \beta}
\sqrt{r(\beta)} \left[F(1)-2\varepsilon_1
\right]\right\},
\end{eqnarray*}
for a constant $\varepsilon_1=\varepsilon_0\sqrt{2 C_0}/16.$ Then,
(\ref{esi16}) and
 (\ref{esi18}) yield $P(\Omega_1)\vee P(\Omega_2)\le 1/4$.
 Moreover,  inequality (\ref{etamber}) implies
\begin{equation*}\begin{array}{l}
\displaystyle {\bf P} \left\{\frac{1}{t \beta^2}\ \log
\left[\sup_{n\le [r(\beta)t]} \sup_{\bx \in P_n}
\int_{S_{n,t}} \beta^{2n}
e^{B_n(\bs ,\bx )
 }  ds_1\cdots ds_n\right]\ge C_0(2-\varepsilon_0/2)\frac{1}{\log(\beta)}\right\}\\[5mm]
 \quad  \ge \displaystyle {\bf P} \left\{\Omega_1^c\cap
 \Omega_2^c\right\} = {\bf P} \left\{(\Omega_1\cup
 \Omega_2)^c\right\} \ge 1 - {\bf P} \left(\Omega_1\right) - {\bf P} \left(
 \Omega_2\right)\ge  \frac12,\end{array},
 \end{equation*}
 which proves (\ref{elower3}), and thus the lower bound of our Theorem \ref{thm:lim-free-energy}.

\subsection{Proof of Theorem \ref{thm:lim-free-energy}: the upper bound}

As in the lower bound section, let us recall that we wish to prove that, for $\beta$ large enough and an arbitrary
positive constant $\varepsilon_0$, we have
$p(\beta)\le C_0 \frac{\beta^2}{\log\left(\beta\right)}\ (1+\varepsilon_0)$.
Again, since $p(\beta)$ exists
and is non-random, we  only need to prove that, for $\beta \ge
\beta_5$ and for $t$ large enough
\begin{equation*}
{\bf P} \left\{\frac{1}{t}\ \log \left[\sum_{n=0}^{\infty} e^{-2 d
t}\sum_{\bx \in P_n} \int_{S_{n,t}} e^{ \beta
B_n(\bs ,\bx )
 }  ds_1\cdots ds_n\right]\le \frac{C_0(1+\varepsilon_0) \beta^2}{\log
 (\beta)}\right\}\ge \frac12.\end{equation*} Actually, we will prove the equivalent inequality
\begin{equation}\label{eupper2}
{\bf P} \left\{\frac{1}{t}\ \log \left[\sum_{n=0}^{\infty} e^{-2 d
t}\sum_{\bx \in P_n} \int_{S_{n,t}} e^{\beta
B_n(\bs ,\bx )
 }  ds_1\cdots ds_n\right]\ge \frac{C_0(1+\varepsilon_0) \beta^2}{\log
 (\beta)}\right\}\le \frac12.\end{equation}

\smallskip

\noindent
{\it Step 1: Setup.}
 Let $\nu=C_0(1+\varepsilon_0) \beta^2/\log
 (\beta)$. For $t$ large enough, the probability defined in (\ref{eupper2}) can be estimated
 from above by the sum of the probabilities of disjoint sets as follows:
\begin{equation}\label{eupper3}
{\bf P} \left\{\sum_{n=0}^{\infty} e^{-2 d
t}\sum_{\bx \in P_n} \int_{S_{n,t}} e^{\beta
B_n(\bs ,\bx )
 }  ds_1\cdots ds_n\ge e^{\nu t}\right\} \le \sum_{l=1}^\infty
 {\bf P} \left(\Lambda_{a_l,b_l}(\nu_l)\right),\end{equation}
 where the sets $\Lambda_{a_l,b_l}(\nu_l)$ are defined by:
 \begin{equation*}
 \Lambda_{a_l,b_l} (\nu_l) = \left\{ \sum_{m=\left[\frac{a_l t \beta^2}{\log^2
(\beta)} +1 \right]}^{\left[\frac{b_l t \beta^2}{\log^2
(\beta)}\right]}  e^{-2d  t} \sum_{\bx \in P_m}
\int_{S_{m,t}} e^{\beta B_m(\bs ,\bx )
 }  ds_1\cdots ds_m\ge e^{\nu_l t}\right\},
 \end{equation*}
and for $l\ge 1$,  the quantities $a_l,b_l,\nu_l$ are of the form
$$a_l=(l-1)\varrho_1, \qquad b_l=l
\varrho_1,\qquad
\nu_l=\left(C_0(1+\varepsilon_0/2)-l\varrho\right)\frac{\beta^2}{\log
 (\beta)},$$
 for two positive constants $\varrho$ and $\varrho_1$ which will be chosen later on,
 at  relation (\ref{defrho}). Notice that
the first set $\Lambda_{a_1,b_1} (\nu_1)$ starts at $m=0$ instead of
$m=\left[a_1 t \beta^2/\log^2 (\beta) +1 \right]$. However, this set
can be handled along the same lines as the other ones. Observe also that the estimate
(\ref{eupper3}) holds true due to the identity
$ \sum_{l=1}^\infty e^{\nu_l t}< e^{\nu t}$, satisfied for $t$ large enough.

\smallskip

\noindent
{\it Step 2: Study the sets} $\Lambda_{a_l,b_l}(\nu_l)$.
Computing the number of terms in each sum and the area of  $S_{m,t}$,
we obtain that
\begin{multline}\label{ear23}
Q_l\equiv\sum_{m=\left[\frac{a_l t \beta^2}{\log^2 (\beta)} +1
\right]}^{\left[\frac{b_l t \beta^2}{\log^2 (\beta)}\right]}
 e^{-2d  t} \sum_{\bx \in P_m}
\int_{S_{m,t}} e^{\beta B_m(\bs ,\bx )
 }  ds_1\cdots ds_m\\
 \le \frac{t(b_l-a_l) \beta^2}{\log^2 (\beta)}\
e^{-2d t} (2d)^{\frac{b_l t \beta^2}{\log^2 (\beta)}}
\exp\left\{\sup_{T_{b_l(\beta)}}
\beta B_m(\bs ,\bx ) \right\}  \sup_{\left[\frac{a_l t
\beta^2}{\log^2 (\beta)} +1 \right] \le m\le \left[\frac{b_l t
\beta^2}{\log^2 (\beta)}\right]} \frac{t^m}{m!},
\end{multline}
where we have set $b_l(\beta)= b_l \beta^2/\log^2
(\beta)$. Furthermore, since $x\le (2d)^x$, we have that
$$
\frac{t(b_l-a_l)\beta^2}{\log^2 (\beta)}\le (2d)^{\frac{b_l t \beta^2}{\log^2
(\beta)}},
$$
and hence the logarithm of $Q_l$ in (\ref{ear23}) is bounded by
\begin{equation}\label{ear24}
\log(Q_l)\le
\frac{2b_l t \beta^2}{\log^2 (\beta)}\log(2d) -2d  t
+\sup_{T_{b_l(\beta)}} \beta B_m(\bs ,\bx )+
\log\sup_{\left[\frac{a_l t \beta^2}{\log^2 (\beta)} +1 \right] \le
m\le \left[\frac{b_l t \beta^2}{\log^2 (\beta)}\right]}
\frac{t^m}{m!}.\end{equation}
Let us find now an estimate for $\sup t^m/m!$ in the expression above:
denote by $K$ a generic constant which depends
only on $d$ and can  change  at each step of our computations. Owing to the bound $(m/3)^m\le m!$,
notice that for $l \ge 2$ and $\beta$
large enough this supremum is attained at the initial point
$\left[\frac{a_l t \beta^2}{\log^2 (\beta)} +1 \right]$. So, we have
\begin{eqnarray*}
\log\sup_{\left[\frac{a_l t \beta^2}{\log^2 (\beta)} +1
\right] \le m\le \left[\frac{b_l t \beta^2}{\log^2 (\beta)}\right]}
\frac{t^m}{m!}
& \le &\frac{a_l t \beta^2}{\log^2 (\beta)} \left[\log t - \log
\left(\frac{a_l t \beta^2}{\log^2 (\beta)}\right) +K \right] \\
& = &\frac{a_l t \beta^2}{\log^2 (\beta)} \left[K - \log (a_l) -
2 \log (\beta) +2\log(\log(\beta)) \right],
\end{eqnarray*}
 for any $l \ge 2$. Using this
fact, the trivial bound $b_l \le 2 a_l$ and plugging the last considerations
 into (\ref{ear24}) and (\ref{ear23}), we end up with
$$
\log(Q_l)\le
\sup_{T_{b_l(\beta)}} \beta
B_m(\bs ,\bx )+ \frac{ a_l t \beta^2}{\log^2
(\beta)} \left[K - 2\log (\beta) - \log (a_l) + 2\log(\log(\beta))
\right].
$$
With this inequality in mind and recalling the definition of $\Lambda_{a_l,b_l}(\nu_l)$, it is readily checked that:
\begin{eqnarray}\label{ear25}
& &{\bf P} \left(\Lambda_{a_l,b_l}(\nu_l)\right)
\\ \nonumber
& &\le  {\bf P}
\left\{\sup_{T_{b_l(\beta)}}
B_m(\bs ,\bx )\ge t\left[\frac{\nu_l}{\beta}
+ \frac{ a_l  \beta}{\log^2 (\beta)} \left[ 2 \log (\beta) + \log
(a_l) -K -2\log(\log(\beta))\right]\right]\right\},
\end{eqnarray}
for any $l \ge 2$.

\smallskip

In order to get an accurate bound on ${\bf P}(\Lambda_{a_l,b_l}(\nu_l))$, we will
now use the information about $B_n$ we have gathered at Section \ref{sec:control-B-n}:
notice for instance that Proposition \ref{pentro} asserts that
\begin{equation*}\lim_{t\rightarrow +\infty} \frac{1}{t}\
\E\sup_{(n,\bs ,\bx )\in T_{b_l(\beta)}} B_n(\bs ,\bx
)=F\left(\frac{b_l \beta^2}{\log^2 (\beta)}\right)=
\frac{\sqrt{b_l} \beta}{\log (\beta)}\ F(1).\end{equation*} Then
Theorem \ref{the22} applied to $\tau^2=t$ together with
(\ref{ear25}) and the last relation imply that, for $l\ge 2$, the
probability of $\Lambda_{a_l,b_l}(\nu_l)$ can be bounded as
\begin{equation}\label{ear28}
{\bf P} \left(\Lambda_{a_l,b_l}(\nu_l)\right)
\le
2 \exp\left\{-\frac{t}{2}
\frac{\beta^2}{\log^2(\beta)} \Xi_l^2\right\},
\end{equation}
where, plugging the definitions of $a_l, b_l$ and $\nu_l$, we can write:
\begin{align*}
&\Xi_l =\frac{\log \beta}{\beta}\left( \frac{\nu_l}{\beta}  +
\frac{a_l \beta}{\log^2 (\beta)} \left[ 2\log (\beta) + \log (a_l)
-K  - 2\log(\log(\beta))\right]
-\frac{\sqrt{b_l} \beta}{\log (\beta)}\ F(1)\right) \\
&= C_0 \left(1 + \frac{\varepsilon_0}{2} \right) -l \varrho
+ \frac{ (l-1) \varrho_1}{\log(\beta)} \left[ 2\log (\beta) + \log
((l-1)\varrho_1) -K - 2\log(\log(\beta))\right] \\
&\hspace{13cm}  - \sqrt{l \varrho_1} F(1)\\
& \ge  C_0
\left(1 + \frac{\varepsilon_0}{2} \right) -l \varrho + 2l
\varrho_1 -2\varrho_1  + \frac{ (l-1) \varrho_1}{\log(\beta)} \left[\log
(\varrho_1) -K - 2\log(\log(\beta))\right]  -  \sqrt{l \varrho_1} F(1).
\end{align*}

\smallskip

Observe that inequality (\ref{ear28}) has been obtained for $l\ge
2$. The same kind of calculations are also valid for $l=1$, except
for the bound on $t^m/m!$. Indeed, in this latter case, owing to
the facts that the maximum of the function $f(x)=\left(
\frac{3t}{x}\right)^x$ is attained at $x=\frac{3t}{e}$ and that
$2d > 3/e$, we obtain that
$$
\log\sup_{0 \le m\le \left[\frac{b_1 t \beta^2}{\log^2
(\beta)}\right]} \frac{t^m}{m!} \le  \frac{3t}{e},
$$
which yield a coefficient $\Xi_1$ of the form:
$$
\Xi_1 := C_0 \left(1 + \frac{\varepsilon_0}{2} \right) - \varrho
+ \varrho_1 \frac{K }{\log(\beta)}- \sqrt{\varrho_1} F(1).
$$
Going back to inequality (\ref{eupper2}), it is also worth mentioning that all
the previous considerations only make sense if $\Xi_l\ge 0$ for all $l\ge 1$.

\smallskip

\noindent
{\it Step 3: Conclusion.}
In order to prove (\ref{eupper2})
and finish the proof of the upper bound, according to (\ref{eupper3})
and (\ref{ear28}), we only need to show, for $t$ large enough, that
\begin{equation}\label{ear28bis}
\sum_{l=1}^\infty  2 \exp\left\{-\frac{t}{2}
\frac{\beta^2}{\log^2(\beta)} \Xi_l^2\right\}<\frac12,\end{equation}
with the additional restriction $\Xi_l\ge 0$. Now, in order to satisfy this latter condition, we  choose $C_0=\frac18\ F(1)^2$, which allows to get rid of the term $-  \sqrt{l \varrho_1} F(1)$ as follows:
\begin{equation*}
\Xi_l  \ge  \left( \sqrt{ C_0 \left( 1 + \frac{\varepsilon_0}{4} \right)} -
\sqrt{ \frac{ 2 l \varrho_1}{1 + \frac{\varepsilon_0}{4}}} \right)^2
+ \left( \frac{ C_0 \varepsilon_0 }{4} - 2 \varrho_1 \right)
 + l \left( 2 \varrho_1 \left( \frac{\varepsilon_0}{4 + \varepsilon_0} \right) - \varrho \right) +  \varrho_1
 \Psi_l,
\end{equation*}
where we have set $\Psi_1:= \frac{K }{\log(\beta)}$, and for any
$l \ge 2$:
$$
\Psi_l:= \frac{ (l-1)}{\log(\beta)} \left[\log
(\varrho_1) -K - 2\log(\log(\beta))\right]. $$
Let us insist at this point on the fact that $C_0=\frac18\ F(1)^2$ is the largest value of $C_0$ allowing such a decomposition.
So, choosing
\begin{equation}\label{defrho}
\varrho_1:=\frac{ C_0 \varepsilon_0}{16},
\quad\mbox{and}\quad
\varrho:= \varrho_1\frac{ \varepsilon_0}{4+ \varepsilon_0},
\end{equation}
we clearly have that
$$
\Xi_l \ge  \frac{ C_0 \varepsilon_0}{8} + l
 \frac{ C_0 \varepsilon_0^2}{16(4+\varepsilon_0)}+  \frac{ C_0 \varepsilon_0}{16}
 \Psi_l.
$$
Thus, there exists $\beta_5$ such that for any $\beta \ge \beta_5$
and for any $l \ge 1$, $\Xi_l$ is strictly positive and the following bound holds true:
\begin{equation*}
\Xi_l \ge  \frac{ C_0 \varepsilon_0}{16} + l  \frac{ C_0 \varepsilon_0^2}{32(4+\varepsilon_0)}.
\end{equation*}
Inequality (\ref{ear28bis}), which ends our proof, follows now easily.

\section{Weak and strong disorder regimes}
\label{sec:weak-strong-disorder}

This section is devoted to the proof of Proposition \ref{prop:weak-strong-disorder},
starting from the result on weak disorder:
\begin{lemma}
Assume $d\ge 3$. Then, for $\beta$ in a neighborhood of 0, the polymer is in the weak disorder regime.
\end{lemma}

\begin{proof}
Similarly to \cite{CY,RT}, it suffices to show that
 \begin{equation}\label{deszeta}
\E[Z_t^2] \le K_1 (\E[Z_t])^2.
\end{equation}
Indeed, recall that the martingale $M$ has been defined by
$M_t= Z_t \exp ( - \frac{\beta^2}{2} t )$. Then, inequality (\ref{deszeta})
yields that $M$ is a martingale bounded in $L^2$ with
$\E[M_t]=1$, and hence $M_\infty=\lim_{t \to \infty} M_t > 0$ on a set of full probability
in $\Omega$, which corresponds to our definition of weak disorder.

\smallskip

In order to check  (\ref{deszeta}), let us consider $b^1$ and $b^2$ two
independent random walks. Then we get
\begin{align*}
\be \left[ Z_t^2 \right] & = E_{b }\left[ \exp\lp \beta ^{2}
\left( t + \int_{0}^{t} \delta_0 (b_s^1  - b_s^2 )
ds \rp \rp \right]\\
&= \left( \be \left[ Z_t \right] \right)^2  E_{b }\left[ \exp\lp
\beta ^{2}  \int_{0}^{t} \delta_0 (b_s^1  - b_s^2 )
ds \rp \right]\\
&\le \left( \be \left[ Z_t \right] \right)^2  E_{ \hat b }\left[
\exp\lp \beta ^{2}  l_\infty (\hat b) \rp \right],
\end{align*}
where $\hat b = b^1 - b^2$ and $l_\infty (\hat b) = \lambda (\{ t
\ge 0, \hat b_t = 0 \})$ with $\lambda$ the Lebesgue measure.
Notice that $\hat b$ is again a continuous time
 random walk  on $\Z^d$ with exponential holding times of
parameter $4d$. Following our notation of Section \ref{sec:intro}, $\hat b$ is described by
its jump times $(\hat \tau_i)_{i \ge 0}$ and its positions $(\hat
x_i)_{i \ge 0}$. Then, if $\beta < \sqrt{4d}$, introducing the obvious
notation $E_{ \hat x }$ and $E_{ \hat \tau }$, we end up with:
\begin{align*}
 &E_{ \hat b }\left[
\exp\lp \beta ^{2}  l_\infty (\hat b) \rp \right]  =
 E_{ \hat b }\left[
\exp\lp \beta ^{2}  \sum_{i=0}^\infty (\hat \tau_{i+1} - \hat \tau_{i}) \delta_0 (\hat x_i )\rp \right]\\
&= E_{ \hat x }\left[  \prod_{i=0}^\infty  E_{ \hat \tau }\left[
 \exp\lp \beta ^{2}  (\hat \tau_{i+1} - \hat \tau_{i}) \delta_0 (\hat x_i )\rp \rp \right]
 =E_{ \hat x }\left[  \prod_{i=0}^\infty  \frac{4d}{4d - \beta^2
\delta_0 (\hat x_i )}   \right]
\\
& =E_{ \hat x }\left[  \left( \frac{4d}{4d - \beta^2 }
\right)^{L_\infty (\hat x)} \right]
 =E_{ \hat x }\left[  \exp \left( \gamma L_\infty (\hat x)
\right) \right],
\end{align*}
where we have set $\gamma=\gamma(\beta)=\log(4d/(4d - \beta^2))$ and
$L_\infty (\hat x) = \#\{ j \le n, \hat x_j=0 \}$, which is the local
time at $x=0$ (and $n=\infty$) of the discrete time random walk induced
by $\hat x$. It is now a well known fact that, in dimension $d\ge 3$ and
for $\gamma$ small enough, we have $E_{ \hat x }\left[  \exp \left( \gamma L_\infty (\hat x)
\right) \right]<\infty$, since $L_\infty (\hat x)$ is a geometric random variable.  Furthermore, $\lim_{\beta\to 0}\gamma(\beta)=0$,
from which our proof is easily finished.

\end{proof}

\smallskip

Let us now say a few words about the remainder of Proposition \ref{prop:weak-strong-disorder}:
point (3) is a direct consequence of our stronger Theorem \ref{thm:lim-free-energy}. As far
as point (2) is concerned, we refer to \cite{CH-dea} for a complete proof of this fact. Like
in \cite{CY,RT}, it is based on an application of It\^o's formula to the medium $W$, which
allows to prove that $\lim_{t\to\infty}\be[M_t^\theta]=0$ for any $\theta\in(0,1)$.

\end{document}